\documentclass[a4paper,11pt]{article}
\usepackage{array}
\usepackage{theorem}
\usepackage{amsmath,amscd,amssymb}
\usepackage{latexsym}
\usepackage{epsfig}

\theorembodyfont{\sl}

\newtheorem{lemma}{Lemma}[section]

\newtheorem{proposition}[lemma]{Proposition}
\newtheorem{theorem}[lemma]{Theorem}

\newtheorem{definition}[lemma]{Definition}

\newcommand{\CC}{\mathbb C}

\newcommand{\HH}{\mathbb H}

\newcommand{\PP}{\mathbb P}

\newcommand{\RR}{\mathbb R}

\newcommand{\ZZ}{\mathbb Z}

\newcommand{\cD}{\mathcal D}

\newcommand{\cF}{\mathcal F}

\newcommand{\cH}{\mathcal H}

\newcommand{\cN}{\mathcal N}

\newcommand{\cS}{\mathcal S}

\newcommand{\To}{\longrightarrow}
\newcommand{\Mapsto}{\mapstochar\longrightarrow}

\renewcommand{\Tilde}{\widetilde}

\renewcommand{\Bar}{\overline}

\newcommand{\SL}{\mathop{\mathrm {SL}}\nolimits}

\newcommand{\Orth}{\mathop{\null\mathrm {O}}\nolimits}

\newcommand{\SOrthtd}{\mathop{\null\mathrm {\widetilde{SO}}}\nolimits}

\renewcommand{\Im}{\mathop{\mathrm {Im}}\nolimits}

\newcommand{\tor}{\mathop{\mathrm {tor}}\nolimits}
\newcommand{\cusp}{\mathop{\mathrm {cusp}}\nolimits}
\newcommand{\vol}{\mathop{\mathrm {vol}}\nolimits}

\newcommand{\divv}{\mathop{\null\mathrm {div}}\nolimits}

\newcommand{\latt}[1]{{\langle{#1}\rangle}}

\newcommand{\Kthree}{\mathop{\mathrm {K3}}\nolimits}
\newcommand{\Ldm}{{L_{2d}^{(m)}}}
\newcommand{\SOrthdm}{\mathop{\null\mathrm {\widetilde{SO}}^+(\Ldm)}\nolimits}
\newcommand{\Fdm}{{\cF}_{2d}^{(m)}}
\newcommand{\SFdm}{{{\cS\cF}_{2d}^{(m)}}}
\newcommand{\Ddm}{{{\cD}_{\Ldm}}}
\newcommand{\Kdm}{{K_{2d}^{(m)}}}
\newcommand{\Ndm}{{N_{2d}^{(m)}}}
\newcommand{\Kmm}{{K_{2}^{(m)}}}
\newcommand{\Tm}{{T_{2,8m+2}}}
\newcommand{\Mm}{{M_{2,8m+2}}}
\newcommand{\vHM}{{\vol_{HM}}}

\newcommand{\qedsymbol}{\mbox{$\Box$}}
\newcommand{\qed}{\unskip\nobreak\hfil\penalty50\hskip1em\hbox{}\nobreak
\hfill\qedsymbol\parfillskip=0pt\finalhyphendemerits=0}
\newenvironment{proof}{\begin{ProofwCaption}{Proof}}{\end{ProofwCaption}}
\newenvironment{ProofwCaption}[1]
 {\addvspace\theorempreskipamount \noindent{\it #1.}\rm}
 {\qed \par \addvspace\theorempostskipamount}

\begin{document}

\title{Hirzebruch-Mumford proportionality and locally symmetric
  varieties of orthogonal type}
\author{V.~Gritsenko, K.~Hulek and G.K.~Sankaran}
\date{}
\maketitle
\begin{abstract}
  For many classical moduli spaces of orthogonal type there are
  results about the Kodaira dimension.  But nothing is known in the
  case of dimension greater than $19$.  In this paper we obtain the
  first results in this direction.  In particular the modular variety
  defined by the orthogonal group of the even unimodular lattice of
  signature $(2,8m+2)$ is of general type if $m\ge 5$.
\end{abstract}

\section{Modular varieties of orthogonal type}
Let $L$ be an integral indefinite lattice of signature $(2,n)$ and
$(\,\,,\,)$ the associated bilinear form.
By $\cD_L$ we denote a connected component of the homogeneous type IV
complex domain of dimension $n$
$$
\cD_L=\{[w]\in\PP(L\otimes\CC)\mid (w,w)=0,\;(w,
\Bar{w})>0\}^+.
$$
$\Orth^+(L)$ is the index $2$ subgroup of the integral orthogonal
group $\Orth(L)$ that leaves $\cD_L$ invariant.  Any subgroup $\Gamma$
of $\Orth^+(L)$ of finite index determines a modular variety
$$
\cF_L(\Gamma)=\Gamma\setminus \cD_L.
$$
By~\cite{BB} this is a quasi-projective variety.

For some special lattices $L$ and subgroups $\Gamma<\Orth^+(L)$ one
obtains in this way the moduli spaces of polarised abelian or Kummer
surfaces ($n=3$, see \cite{GH}), the moduli space of
Enriques surfaces ($n=10$, see \cite{BHPV}), and the moduli
spaces of polarised or lattice-polarised $\Kthree$ surfaces 
($0< n\le 19$, see \cite{Nik1, Dol}).  Other interesting modular
varieties of orthogonal type include the period domains of the
symplectic manifolds and certain varieties associated to fermionic and
bosonic strings.

It is natural to ask about the birational type of $\cF_L(\Gamma)$. For
many classical moduli spaces of orthogonal type there are results
about the Kodaira dimension, but nothing is known in the case of
dimension greater than $19$.  In this paper we obtain the first
results in this direction.  We determine the Kodaira dimension of many
quasi-projective varieties associated with two series of even
lattices.

To explain what these varieties are, we first introduce the {\it
  stable orthogonal group\/} $\Tilde{\Orth}(L)$ of a nondegenerate
even lattice $L$. This is defined (see \cite{Nik2} for more details) 
to be the subgroup of $\Orth(L)$
which acts trivially on the discriminant group $A_L=L^\vee/L$, where
$L^\vee$ is the dual lattice. If $\Gamma<\Orth(L)$ then we write
$\Tilde\Gamma=\Gamma\cap \Tilde\Orth(L)$. Note that if $L$ is
unimodular then $\Tilde\Orth(L)=\Orth(L)$. 

The first series of varieties we want to study, which we call the
modular varieties of {\it unimodular type}, is
\begin{equation}\label{mvar}
{\cF}_{II}^{(m)} = {\Orth}^+(II_{2,8m+2})
\backslash \cD_{II_{2,8m+2}}.
\end{equation}
$\cF_{II}^{(m)}$ is of dimension $8m+2$ and arises from the even
unimodular lattice of signature $(2, 8m+2)$
\begin{equation*}
II_{2,8m+2}=2U\oplus mE_8(-1),
\end{equation*}
where $U$ denotes the hyperbolic plane and $E_8(-1)$ is the negative
definite lattice associated to the root system $E_8$.
The case $m=3$ is of
particular interest: $\cF_{II}^{(3)}$ is of dimension~$26$ and arises
in the context of bosonic strings.

The second series, which we call the modular varieties of 
{\it K3 type}, is
\begin{equation}\label{2dvar}
{\cF}_{2d}^{(m)} = \widetilde{\Orth}^+(L_{2d}^{(m)})
\backslash \cD_{L_{2d}^{(m)}}.
\end{equation}
${\cF}_{2d}^{(m)}$ is of dimension $8m+3$ and arises from the lattice
\begin{equation*}
L_{2d}^{(m)}=2U\oplus mE_8(-1)\oplus \latt{-2d},
\end{equation*}
where $\latt{-2d}$
denotes a lattice generated by a vector of square $-2d$.

The first three members of the series $\cF_{2d}^{(m)}$ have
interpretations as moduli spaces. $\cF_{2d}^{(2)}$ is the moduli space
of polarised $K3$ surfaces of degree $2d$.  For $m=1$ the
$11$-dimensional variety $\cF_{2d}^{(1)}$ is the moduli space of
lattice-polarised $\Kthree$ surfaces, where the polarisation is
defined by the hyperbolic lattice $\latt{2d}\oplus E_8(-1)$ (see
\cite{Nik1, Dol}).  For $m=0$ and $d$ prime the $3$-fold $\cF_{2d}^{(0)}$ is
the moduli space of polarised Kummer surfaces (see \cite{GH}).

\begin{theorem}\label{gentype}
The modular varieties of unimodular and $\Kthree$ type are varieties of
general type if $m$ and $d$ are sufficiently large. More precisely:
\begin{itemize}
\item[(i)]
If $m\ge 5$ then  the modular varieties ${\cF}_{II}^{(m)}$
and ${\cF}_{2d}^{(m)}$ (for any $d\ge 1$) are of general type.
\item[(ii)]
For $m=4$ the varieties  ${\cF}_{2d}^{(4)}$ are of general type if $d\ge 3$
and $d\ne 4$.
\item[(iii)]
For $m=3$ the varieties ${\cF}_{2d}^{(3)}$ are of general type if  $d\ge 1346$.
\item[(iv)]
For $m=1$ the varieties ${\cF}_{2d}^{(1)}$ are of general type if  $d\ge 1537488$.
\end{itemize}
\end{theorem}
{\bf Remark.}  The methods of this paper are also applicable if $d=2$.
Using them, one can show that the moduli space  ${\cF}_{2d}^{(2)}$ of polarised $\Kthree$
surfaces of degree $2d$
is of general type if $d\ge
231000$.  This case was studied in \cite{GHS2}, where, using a
different method involving special pull-backs of the Borcherds
automorphic form $\Phi_{12}$ on the domain $\cD_{II_{2,26}}$, we
proved that ${\cF}_{2d}^{(2)}$ is of general type if $d>61$ or $d=46$,
$50$, $54$, $57$, $58$, $60$.

The methods of~\cite{GHS2} do not appear to be applicable in the other
cases studied here. Instead, the proof of Theorem~\ref{gentype}
depends on the existence of a good toroidal compactification of
$\cF_L(\Gamma)$, which was proved in \cite{GHS2}, and on the exact
formula for the Hirzebruch-Mumford volume of the orthogonal group
found in \cite{GHS1}.

We shall construct pluricanonical forms on a suitable compactification
of the modular variety $\cF_{L}(\Gamma)$ by means of modular forms.  Let
$\Gamma < \Orth^+(L)$ be a subgroup of finite index, which naturally
acts on the affine cone $\cD_{L}^{\bullet}$ over $\cD_{L}$. In what follows
we assume that $\dim\cD_{L} \ge 3$.

\begin{definition}\label{modforms}
  A modular form of weight $k$ and character $\chi\colon \Gamma\to
  \CC^*$ with respect to the group $\Gamma$ is a holomorphic function
$$
F: \cD_{L}^{\bullet} \to \CC
$$
which has the two properties
$$
\begin{array}{rcll}
F(tz)&=&t^{-k}F(z)& \forall\ t\in \CC^*,\\
F(g(z)) & = &\chi(g)F(z) & \forall\ g \in \Gamma.
\end{array}
$$
\end{definition}
The space of modular forms is denoted by $M_k(\Gamma,\chi)$. The
space of cusp forms, {i.e.} modular forms vanishing on the
boundary of the Baily--Borel compactification of $\Gamma\backslash
\cD_L$, is denoted by $S_k(\Gamma,\chi)$. We can reformulate the
definition of modular forms in geometric terms.  Let $F\in
M_{kn}(\Gamma,\det^k)$ be a modular form, where $n$ is the dimension
of $\cD_{L}$. Then
$$
F(dZ)^k\in H^0(\cF_{L}(\Gamma)^\circ, \Omega^{\otimes k}),
$$
where $dZ$ is a holomorphic volume form on $\cD_{L}$,
$\Omega$ is the sheaf of germs of canonical  $n$-forms
on $\cF_{L}(\Gamma)$ and $\cF_{L}(\Gamma)^\circ$ is the open smooth part
of $\cF_{L}(\Gamma)$ such that the projection
$
\pi:\cD_{L}\to  \Gamma\backslash \cD_{L}
$
is unramified over $\cF_{L}(\Gamma)^\circ$.

The main question in the proof of Theorem~\ref{gentype} is how to
extend the form $F(dZ)^k$ to $\cF_{L}(\Gamma)$ and to a suitable
toroidal compactification~$\cF_L(\Gamma)^{\tor}$.
There are three possible kinds of obstruction to this, which we call
(as in~\cite{GHS2}) elliptic, reflective and cusp
obstructions. Elliptic obstructions arise if
$\cF_{L}(\Gamma)^{\tor}$ has
non-canonical singularities arising from fixed loci of the action
of the group $\Gamma$.  Reflective obstructions arise because the projection $\pi$ is
branched along divisors whose general point is smooth in
$\cF_{L}(\Gamma)$. Cusp obstructions arise when we extend the form
from $\cF_L(\Gamma)$ to $\cF_L(\Gamma)^{\tor}$.

The problem of elliptic obstructions was solved for $n\ge 9$ in \cite{GHS2}.
\begin{theorem}\label{torcomp}(\cite[Theorem 2.1]{GHS2})
  Let $L$ be a lattice of signature $(2,n)$ with $n\ge 9$, and let
  $\Gamma<\Orth^+(L)$ be a subgroup of finite index. Then there exists
  a toroidal compactification $\cF_L(\Gamma)^{\tor}$ of
  $\cF_L(\Gamma)=\Gamma\backslash\cD_L$ such that $\cF_L(\Gamma)^{\tor}$
  has canonical singularities and there are no branch divisors in the
  boundary.  The branch divisors in $\cF_L(\Gamma)$ arise from the
  fixed divisors of reflections.
\end{theorem}
Reflective obstructions, that is branch divisors, are a special
problem related to the orthogonal group. They do not appear in the
case of moduli spaces of polarised abelian varieties of dimension
greater than $2$, where the modular group is the symplectic group.
There are no quasi-reflections in the symplectic group even for $g=3$.

The branch divisor is defined by special reflective vectors
in the lattice $L$.  This description is given in \S 2.
To estimate the reflective obstructions
we use the Hirzebruch-Mumford proportionality principle
and  the exact formula for the  Hirzebruch-Mumford volume
of the orthogonal group found in \cite{GHS1}.
We do the numerical estimation in \S 4.

We treat the cusp obstructions in \S 3, using special cusp forms of low weight
(the lifting of Jacobi forms) constructed in \cite{G2} 
and the low-weight cusp form trick (see \cite{G2} and \cite{GHS2}).

\section{The branch divisors}

To estimate the obstruction to extending pluricanonical
forms to a smooth projective model of $\cF_L(\Tilde\Orth^+(L))$
we have to determine the branch divisors of the projection
\begin{equation}\label{pi}
\pi\colon \cD_L\to \cF_L(\Tilde\Orth^+(L))=\Tilde{\Orth}^+(L)\setminus \cD_L.
\end{equation}
According to \cite[Corollary 2.13]{GHS2} these divisors
are defined by reflections $\sigma_r\in \Orth^+(L)$, where
$$
\sigma_r(l)=l-\frac{2(l,r)}{(r,r)}r,
$$
coming from vectors $r\in L$ with $r^2<0$ that are {\it stably
  reflective\/}: by this we mean that $r$ is primitive and
$\sigma_r$ or $-\sigma_r$ is in $\Tilde\Orth^+(L)$. 
By a {\it $(k)$-vector\/} for $k\in\ZZ$ we mean a primitive vector $r$ with $r^2=k$.

Let $D$ be the exponent of the finite abelian group $A_L$
and let the {\it divisor\/} $\divv(r)$ of $r\in L$ be the positive generator
of the ideal $(l,L)$.  We note that $r^*=r/\divv(r)$ is a primitive
vector in $L^\vee$.  In \cite[Propositions 3.1--3.2]{GHS2} we proved
the following.
\begin{lemma}\label{rsigma}\textbf{}
Let $L$ be an even integral lattice
of signature $(2,n)$.
If $\sigma_r\in \Tilde\Orth^+(L)$ then
$r^2=-2$.
If $-\sigma_r\in \Tilde\Orth^+(L)$, then
$r^2=-2D$ and $\divv(r)=D\equiv 1$ mod~$2$
or $r^2=-D$ and $\divv(r)=D$ or $D/2$.
\end{lemma}

We need also the following well-known property of the stable
orthogonal group.
\begin{lemma}\label{stable}
  For any sublattice $M$ of an even lattice $L$ the group $\Tilde\Orth(M)$
  can be considered as a subgroup of $\Tilde\Orth(L)$.
\end{lemma}
\begin{proof}
Let $M^\perp$ be the orthogonal complement of $M$ in $L$.
We have as usual
$$
M\oplus M^\perp \subset L\subset L^{\vee}\subset M^{\vee}\oplus (M^\perp)^{\vee}.
$$
We can extend $g\in \Tilde\Orth(M)$ to $M\oplus M^\perp$ by putting
$g|_{M^\perp}\equiv \hbox{id}$. It is clear that $g\in
\Tilde\Orth(M\oplus M^\perp)$. 
For any $l^\vee\in L^{\vee}$ we have $g(l^\vee)\in l^\vee+(M\oplus
M^\perp)$.  In particular, $g(l)\in L$ for any $l\in L$ and $g\in
\Tilde\Orth(L)$.
\end{proof}

We can describe the components of the branch locus in terms of
homogeneous domains. For $r$ a stably reflective vector in $L$ we put
$$
H_r=\{[w]\in \PP(L\otimes \CC)\mid (w,r)=0\},
$$
and let $\cN$ be the union of all hyperplane sections
$H_r\cap \cD_L$ over all stably reflective vectors~$r$.
\begin{proposition}\label{branchdiv}
  Let $r\in L$ be a stably reflective vector: suppose that $r$ and $L$ do
  not satisfy $D=4$, $r^2=-4$, $\divv(r)=2$.  Let $K_r$ be
  the orthogonal complement of $r$ in $L$. Then the associated 
component $\pi(H_r\cap\cD_L)$ of the
  branch locus $\cN$ is of the form
  ${\Tilde{\Orth}^+(K_{r})}\backslash \cD_{K_r}$.
\end{proposition}
\begin{proof}
We have $H_r\cap \cD_L=\PP(K_r)\cap \cD_L=\cD_{K_r}$. Let
\begin{equation}\label{GammaK}
\Gamma_{K_r}=\{ \varphi \in \Tilde{\Orth}^+(L)\mid \varphi(K_r)=K_r\}.
\end{equation}
$\Gamma_{K_r}$ maps to a subgroup of ${\Orth}^+(K_r)$. 
The inclusion of $\Tilde\Orth(K_r)$ in $\Tilde\Orth(L)$ (Lemma~\ref{stable}) 
preserves the spinor norm (see~\cite[\S 3.1]{GHS1}),  because $K_r$ has signature
$(2,n-1)$ and so $\Tilde{\Orth}^+(K_r)$ becomes a
subgroup of $\Tilde\Orth^+(L)$.

Therefore the image of $\Gamma_{K_r}$ contains
$\Tilde\Orth^+(K_r)$ for any $r$. Now we prove that this image
coincides with $\Tilde\Orth^+(K_r)$ for all $r$, except perhaps if
$D=4$, $r^2=-4$ and $\divv(r)=2$.

Let us consider the inclusions
$$
\latt{r}\oplus K_r\subset L\subset L^\vee
\subset \latt{r}^\vee\oplus K_r^{\vee}.
$$
By standard arguments (see \cite[Proposition 3.6]{GHS2}) we see that
$$
|\det K_r|= \frac{|\det L|\cdot |r^2|}{\divv(r)^2}\quad
\text{and}\quad
[L:\latt{r} \oplus K_r]=\frac{|r^2|}{\divv(r)}=1\ \text{ or }\ 2.
$$
If the index is $1$, then it is clear that the image of
$\Gamma_{K_r}$ is $\Tilde\Orth^+(K_r)$.
Let us assume that the index is equal to $2$.
In this case the lattice $\latt{r}^\vee$ is generated by
$r^\vee=-r/(r,r)=r^*/2$, where $r^*=r/\divv(r)$ is a primitive
vector in $L^\vee$.
In particular $r^\vee$ represents a non-trivial class
in $\latt{r}^\vee\oplus K_r^{\vee}$ modulo $L^\vee$.
Let us take $k^\vee \in K_r^\vee$ such that $k^\vee\not\in L^\vee$.
Then $k^\vee+r^\vee\in L^\vee$ and
$$
\varphi(k^\vee)-k^\vee\equiv r^\vee-\varphi(r^\vee)\mod L.
$$
We note that if $\varphi \in\Gamma_{K_r}$ then
$\varphi(r)=\pm r$. Hence
$$
\varphi(k^\vee)-k^\vee\equiv
\begin{cases}
0\mod L&\text{if } \varphi(r)=r\\
r^*\mod L&\text{if } \varphi(r)=-r
\end{cases}
$$
Since $\varphi(r^*)\equiv r^*\mod L$, we cannot have
$\varphi(r)=-r$ unless $\divv(r)=1$ or $2$.
Therefore we have proved that
$\varphi(k^\vee)\equiv k^\vee\mod K_r$
($K_r=K_r^\vee\cap L$), except possibly if $D=4$, $r^2=-4$, $\divv(r)=2$.
\end{proof}
The group $\Tilde{\Orth}^+(L)$ acts on $\cN$. We need to estimate the
number of components of $\Tilde{\Orth}^+(L)\setminus \cN$.  This will
enable us to estimate the reflective obstructions to extending
pluricanonical forms which arise from these branch loci.

For the even unimodular lattice $II_{2,8m+2}$  any  primitive vector
$r$ has $\divv(r)=1$. Consequently $r$ is stably reflective 
if and only if $r^2=-2$.

For $L_{2d}^{(m)}$ the reflections and the corresponding branch
divisors arise in two different ways, according to Lemma~\ref{rsigma}.
We shall classify the orbits of such vectors.

\begin{proposition}\label{orbits}
Suppose $d$ is a positive integer.
\begin{itemize}
\item[\rm{(i)}]
Any two $(-2)$-vectors in the lattice $II_{2,8m+2}$ are
equivalent modulo ${\Orth}^+(II_{2,8m+2})$, and
the orthogonal complement of a $(-2)$-vector $r$ is isometric to
$$
K_{II}^{(m)}=U\oplus mE_8(-1)\oplus \latt{2}.
$$
\item[\rm{(ii)}]
There is one $\Tilde{\Orth}^+(L_{2d}^{(m)})$-orbit of $(-2)$-vectors $r$
in $L_{2d}^{(m)}$ with $\divv(r)=1$. If $d\equiv 1$ mod~$4$ 
then there is a second orbit of $(-2)$-vectors, with $\divv(r)=2$.
The orthogonal complement of a $(-2)$-vector $r$ in $L_{2d}^{(m)}$
is isometric to
$$
K_{2d}^{(m)}=  U\oplus mE_8(-1)\oplus \latt{2} \oplus\latt{-2d} ,
$$
if $\divv(r)=1$, and to
$$
N_{2d}^{(m)}=U\oplus m E_8 (-1) \oplus
\begin{pmatrix}
1& 2\\
\frac{1-d}{2} &  1
\end{pmatrix},
$$
if $\divv(r)=2$.
\item[\rm{(iii)}] 
The orthogonal complement of a $(-2d)$-vector $r$ in $L_{2d}^{(m)}$
is isometric to
$$
II_{2,8m+2}= 2U\oplus mE_8(-1)
$$ 
if $\divv(r)=2d$, and to
$$
K_{2}^{(m)} = U\oplus mE_8(-1)\oplus \latt{2}\oplus \latt{-2}
\quad\text{or}\quad
T_{2,8m+2}=U\oplus U(2)\oplus  mE_8(-1)
$$
if $\divv(r)=d$.
\item[\rm{(iv)}] Suppose $d>1$.  The number of
  $\Tilde{\Orth}(L_{2d}^{(m)})$-orbits of $(-2d)$-vectors
  with $\divv(r)=2d$ is $2^{\rho(d)}$.  The number of
  $\Tilde{\Orth}(L_{2d}^{(m)})$-orbits of $(-2d)$-vectors
  with $\divv(r)=d$ is
$$
\begin{cases}
2^{\rho(d)}&    \text{if $d$ is odd or }\  d\equiv 4 \mod 8;\\
2^{\rho(d)+1}&  \text{if }\  d\equiv 0 \mod 8;\\
2^{\rho(d)-1}&   \text{if }\  d\equiv 2 \mod 4.
\end{cases}
$$
\end{itemize}
Here $\rho(d)$ is the number of prime divisors of $d$.
\end{proposition}
\begin{proof}
If the lattice $L$ contains two hyperbolic planes
then according to the well-known result
of Eichler (see \cite[\S 10]{E}) the $\Tilde{\Orth}^+(L)$-orbit
of a primitive vector $l\in L$ is completely defined by two invariants:
by its length $(l,l)$ and by its image $l^*+L$ in the discriminant group $A_L$,
where $l^*=l/\divv(l)$.

i) If $u$ is a primitive vector of an even unimodular lattice $II_{2,8m+2}$
then $\divv(u)=1$ and there is only one $\Orth(II_{2,8m+2})$-orbit
of $(-2)$-vectors. Therefore we can take $r$ to be a $(-2)$-vector in
$U$, and the form of the orthogonal complement is obvious.

ii) In the lattice $L_{2d}^{(m)}$
we  fix a generator $h$ of its $\langle -2d\rangle$-part.
Then for any $r\in L_{2d}^{(m)}$ we can write $r=u+xh$,
where $u\in II_{2,8m+2}$ and $x\in \ZZ$.
It is clear that $\divv(r)$ divides $r^2$.
If $f|\divv(r)$, where $f=2$, $d$ or $2d$,
then the vector $u$ is also divisible by $f$.
Therefore the $(-2)$-vectors form two possible orbits of vectors
with divisor equal to $1$ or $2$.
If $r^2=-2$ and $\divv(r)=2$ then $u=2u_0$ with $u_0\in 2U\oplus mE_8(-1)$
and we see that in this case  $d\equiv 1\mod 4$.
This gives us two different orbits for such $d$.
In both cases we can find a $(-2)$-vector $r$ in the sublattice
$U\oplus \langle -2d\rangle$.
Elementary calculation gives us the orthogonal complement of $r$.

iii) This was proved in \cite[Proposition 3.6]{GHS2}
for $m=2$. For general $m$ the proof is the same.

iv) To find the number of orbits of $(-2d)$-vectors
we have to consider two cases.

a) Let $\divv(r)=2d$. Then $r=2du+xh$ and $r^*\equiv (x/2d)h \mod L$,
where $u\in II_{2,8m+2}$ and $x$ is modulo $2d$.  Moreover
$(r,r)=4d^2(u,u)-x^22d=-2d$. Thus $x^2\equiv 1 \mod 4d$.  This
congruence has $2^{\rho(d)}$ solutions modulo $2d$.  For any such $x
\mod 2d$ we can find a vector $u$ in $2U\oplus mE_8(-1)$ with
$(u,u)=(x^2-1)/2d$.  Then $r=2du+xh$ is primitive (because $u$ is not
divisible by any divisor of $x$) and $(r,r)=-2d$.

b) Let $\divv(r)=d$.  Then $r=du+xh$, where $u$ is primitive, $r^*
\equiv (x/d)h \mod L$ and $x$ is modulo $d$.  We have $(r^*,r^*)\equiv
-2x^2/d \mod 2\ZZ$ and $x^2\equiv 1 \mod d$.  For any solution modulo
$d$ we can find as above $u\in 2U\oplus mE_8(-1)$ such that $r=du+xh$ is
primitive and $(r,r)=-2d$.  It is easy to see that the number of solutions $\{x \mod
d\, |\, x^2\equiv 1 \mod d\}$ is as stated.
\end{proof}
{\bf Remark}. To calculate the number of the
branch divisors arising from vectors $r$ with $r^2=-2d$ one has to
divide the corresponding number of orbits found in
Proposition~\ref{orbits}(iv) by $2$ if $d>2$.  This is because $\pm r$
determine different orbits but the same branch divisor.
For $d=2$ the proof shows that there is one divisor
for each orbit given in Proposition~\ref{orbits}(iv).

\section{Modular forms of low weight}
In this section we let $L=2U\oplus L_0$ be an even lattice of signature
$(2,n)$ with two hyperbolic planes.  We choose a primitive isotropic
vector $c_1$ in $L$.  This vector determines a $0$-dimensional cusp
and a tube realisation of the domain $\cD_L$.  The tube domain (see
$\cH(L_1)$ below) is a complexification of the positive cone of the
hyperbolic lattice $L_1=c_1^{\perp}/c_1$.  If $\divv(c_1)=1$ we call this cusp
{\it standard\/} (as above, by~\cite{E} there is only one standard cusp). 
In this case  $L_1=U\oplus L_0$.
In \cite[\S 4]{GHS2} we proved that any $1$-dimensional boundary component
of $\Tilde{\Orth}^+(L)\setminus \cD_L$ contains the standard $0$-dimensional cusp
if every isotropic (with respect to the discriminant form: see \cite [\S{}1.3]{Nik2})
subgroup of $A_L$ is cyclic.

Let us fix a $1$-dimensional cusp by choosing two copies of $U$ in
$L$.  (One has to add to $c_1$ a primitive isotropic vector $c_2\in
L_1$ with $\divv(c_2)=1$).  Then $L=U\oplus L_1=U\oplus(U\oplus L_0)$
and the construction of the tube domain
may be written down simply in coordinates. We have
$$
\cH(L_1)=\cH_{n}=\{Z=(z_n, \dots, z_1)
\in \HH_1\times \CC^{n-2}\times\HH_1; \; (\Im Z,\Im Z)_{L_1}>0\},
$$
where $Z\in L_1\otimes \CC$ and $(z_{n-1},\dots, z_2)\in L_0\otimes \CC$.
(We represent $Z$ as a column vector.)
An isomorphism between $\cH_{n}$
and $\cD_L$ is given by
\begin{eqnarray}\label{tubedomain}
p\colon \cH_{n}&\To& \cD_L\\
Z=(z_n, \dots, z_1)
&\Mapsto& \big( -\frac{1}{2}(Z, Z)_{L_1}:
z_n: \cdots :z_1 : 1\big)\nonumber.
\end{eqnarray}
The action of $\Orth^+(L\otimes \RR)$ on $\cH_{n}$ is given
by the usual fractional linear transformations.
A calculation shows that the Jacobian of the transformation of $\cH_{n}$ defined
by $g\in \Orth^+(L\otimes \RR)$ is equal to $\det(g)j(g,Z)^{-n}$,
where $j(g,Z)$ is the last ($(n+2)$-nd) coordinate of
$g\big(p(Z)\big)\in\cD_L$.  Using this we define the automorphic
factor
\begin{eqnarray*}
J\colon  \Orth^+(L\otimes \RR)\times \cH_{n+2}&\to& \CC^*\\
(g, Z)&\mapsto& (\det g)^{-1} \cdot j (g, Z)^n.
\end{eqnarray*}
The connection with pluricanonical forms is the following. 
Consider the form
$$
dZ=dz_1\wedge\cdots\wedge d z_{n}\in
\Omega^{n}(\cH_{n}).
$$
$F(dZ)^k$ is a $\Gamma$-invariant $k$-fold pluricanonical form on
$\cH_{n}$, for $\Gamma$ a subgroup of finite index of $\Orth^+(L)$, if
$F\big(g(Z)\big)=J(g, Z)^{k}F(Z)$ for any $g\in \Gamma$; in other
words if $F\in M_{nk}(\Gamma, \det^k)$ (see Definition~\ref{modforms}).  To prove
Theorem~\ref{gentype} we need cusp forms of weight smaller than the dimension
of the corresponding modular variety.
\begin{proposition}\label{cuspforms}
For unimodular type, cusp forms of weight $12+4m$ exist: that is
$$
\dim S_{12+4m}(\Orth^+(II_{2,8m+2}))>0.
$$
For $\Kthree$ type we have the bounds
\begin{eqnarray*}
\dim S_{11+4m}(\Tilde{\Orth}^+(L_{2d}^{(m)}))&>&0 \ \text{ if }\ d>1;\\
\dim S_{10+4m}(\Tilde{\Orth}^+(L_{2d}^{(m)}))&>&0 \ \text{ if }\ d\ge 1;\\
\dim S_{7+4m}(\Tilde{\Orth}^+(L_{2d}^{(m)}))&>&0\ \text{ if }\ d\ge 4;\\
\dim S_{6+4m}(\Tilde{\Orth}^+(L_{2d}^{(m)}))&>&0\ \text{ if }\ d=3\
\text{ or }\ d\geq 5;\\
\dim S_{5+4m}(\Tilde{\Orth}^+(L_{2d}^{(m)}))&>&0\ \text{ if }\ d=5\ \text{ or }\ d\ge 7;\\
\dim S_{2+4m}(\Tilde{\Orth}^+(L_{2d}^{(m)}))&>&0\ \text{ if }\
d>180.
\end{eqnarray*}
\end{proposition}
\begin{proof}
For any $F(Z)\in M_k(\Tilde{\Orth}^+(L))$ we can consider its
Fourier-Jacobi expansion at the  $1$-dimensional cusp
fixed above
$$
F(Z)=f_0(z_1)+\sum_{m\ge 1}f_m(z_1; z_2,\dots z_{n-1})
\exp(2\pi i m z_n).
$$
A lifting construction of modular
forms $F(Z)\in M_k(\Tilde{\Orth}^+(L))$ with trivial character by
means of the first Fourier--Jacobi coefficient is given in \cite{G1}, \cite{G2}.  We note that
$f_1(z_1; z_2,\dots, z_{n-1})\in J_{k,1}(L_0)$, where $J_{k,1}(L_0)$ is
the space of the Jacobi forms of weight $k$ and index $1$.  A more
general construction of the additive lifting was given in \cite{B2} but
for our purpose the construction of \cite{G2} is sufficient.

The dimension of $J_{k,1}(L_0)$ depends only on
the discriminant form and  the rank of $L_0$ (see \cite[Lemma 2.4]
{G2}).  In particular, for the special cases of $L=II_{2,8m+2}$ and
$L=L_{2d}^{(m)}$ we have
$$
J_{k+4m,1}^{\cusp}(mE_8(-1))\cong S_{k}(\SL_2(\ZZ))
$$
and
$$
J_{k+4m,1}^{\cusp}(mE_8(-1)\oplus\langle-2d\rangle)\cong
J_{k,d}^{\cusp},
$$
where $J_{k,d}^{\cusp}$ is the space of the usual Jacobi cusp forms
in two variables of weight $k$ and index $d$ (see [EZ]) and $
S_{k}(\SL_2(\ZZ))$ is the space of weight $k$ cusp forms for $\SL_2(\ZZ)$.

The lifting of a Jacobi cusp form of index one is a cusp form of
the same weight with respect to $\Orth^+(II_{2,8m+2})$ or
$\Tilde{\Orth}^+(L_{2d}^{(m)})$ with trivial character.  The fact that
we get a cusp form was proved in \cite{G2} for maximal lattices, i.e., if
$d$ is square-free.  In \cite[\S 4]{GHS2} we extended this to
all lattices $L$ for which the isotropic subgroups of the discriminant
$A_L$ are all cyclic, which is true in all cases considered here.

To prove the unimodular type case of Proposition~\ref{cuspforms} we can take the Jacobi
form corresponding to the cusp form $\Delta_{12}(\tau)$. Using the
Jacobi lifting construction we obtain a cusp form of weight $12+4m$ with
respect to $\Orth^+(II_{2,8m+2})$.

For the $\Kthree$ type case we need the dimension formula for the space of Jacobi cusp forms
$J_{k,d}^{\cusp}$ (see~\cite{EZ}). For a positive integer $l$ one sets
$$
\{l\}_{12}=
\left\{\begin{array}{ll}
 \lfloor\frac{l}{12}\rfloor &\mbox{ if } l
\not\equiv 2\bmod 12\\
\\
\lfloor\frac{l}{12}\rfloor -1
&\mbox{ if } l\equiv 2\bmod 12.
\end{array} \right.
$$
Then if $k>2$ is even
$$
\dim  J_{k,d}^{\,\cusp}
=\sum\limits^{d}_{j=0}\left(\{k+2j\}_{12}-\left\lfloor\frac{j^2}{4d}
\right\rfloor\right),
$$
and if $k$ is odd
$$
\dim  J_{k,d}^{\,\cusp}
=\sum\limits^{d-1}_{j=1}\left(\{k-1+2j\}_{12}-\left\lfloor\frac{j^2}{4d}
\right\rfloor\right).
$$
This gives the bounds claimed. For $k=2$, using the results
of~\cite{SZ} one can also calculate $\dim
J_{2,d}^{\cusp}$: there is an extra term, $\lceil \sigma_0(d)/2\rceil$,
where $\sigma_0(d)$ denotes the number of divisors of $d$. 
This gives
$\dim J_{2,d}^{\cusp}>0$ if $d>180$ and for some smaller values of~$d$. 
\end{proof}

\section{Kodaira dimension results}

In this section we prove Theorem~\ref{gentype}.  Our strategy is the
following.  For $\Gamma\subseteq \Tilde\Orth^+(L)$ we choose a cusp 
form $F_a\in S_a(\Gamma)$ of low weight $a$, 
i.e. $a$ strictly less than the dimension $n$ of
$\cF_L(\Gamma)$.
Then we consider elements $F\in F^k_{a}M_{k(n-a)}(\Gamma, \det^k)$:
for simplicity we assume that $k$ is even. Such an $F$ vanishes to
order at least $k$ on the boundary of any toroidal compactification.
Hence if $dZ$ is the volume element on $\cD_L$ defined in \S 3
it follows that $F(dZ)^k$ extends as a $k$-fold pluricanonical form to
the general point of every boundary component of $\cF_L(\Gamma)^{\tor}$.
Now assume that we have chosen the toroidal compactification so that
all singularities are canonical and that there is no ramification
divisor which is contained in the boundary. Such toroidal
compactifications exist by Theorem~\ref{torcomp} if the dimension
$n\ge 9$. Then the only obstructions to extending $F(dZ)^k$ to a
smooth projective model are the reflective obstructions, coming from
the ramification divisor of the quotient map $\pi\colon \cD_{L}\to
\cF_L(\Gamma)$ studied in \S 2.

Let $\cD_K$ be an irreducible component of this ramification divisor.
Recall from Proposition~\ref{branchdiv} that
$\cD_K=\PP(K\otimes\CC)\cap\cD_L$ where $K=K_r$ is the
orthogonal complement of a stably reflective vector $r$.  For the lattices
chosen in Theorem~\ref{gentype} all irreducible components of the
ramification divisor are given in Proposition~\ref{orbits}.

\begin{proposition}\label{obstructions}
We assume that $k$ is even and that the dimension $n\ge 9$. 
For $\Gamma\subseteq \Tilde\Orth^+(L)$, 
the obstruction to extending forms
$F(dZ)^k$ where $F\in F^k_{a}M_{k(n-a)}(\Gamma)$ to
$\cF_L(\Gamma)^{\tor}$ lies in the space
$$
B=\bigoplus_{K}B(K)
=\bigoplus_{K}\bigoplus_{\nu = 0}^{k/2-1} M_{k(n-a)+2\nu}({\Gamma\cap\Tilde{\Orth}^+(K)}),
$$
where the direct sum is taken over all
irreducible components $\cD_K$ of the ramification divisor of the
quotient map $\pi\colon \cD_{L}\to \cF(\Gamma)$.
\end{proposition}
\begin{proof}
Let $\sigma \in \Gamma$ be plus or minus a reflection
whose fixed point locus is $\cD_K$.
We can extend the differential form provided that $F$
vanishes of order $k$ along every irreducible component $\cD_K$ of the
ramification divisor.

If $F_{a}$ vanishes along $\cD_K$ then $K$ gives no restriction on the
second factor of the modular form $F$.

Now let $\{w=0\}$ be a local equation for $\cD_K$. Then
$\sigma^*(w)=-w$ (this is independent of whether
$\sigma$ or $-\sigma$ is the reflection).
For every modular  form $F\in M_{k}(\Gamma)$ of even weight we have
$F\big(\sigma(z)\big)=F(z)$.
This implies that if $F(z)\equiv 0$ on $\cD_K$,
then $F$ vanishes to even order on $\cD_K$.

We denote by $M_{2b}(\Gamma)(-\nu\cD_K)$ the space of modular forms of
weight $2b$ which vanish of order at least $\nu$ along $\cD_K$.  Since
the weight is even we have $M_{2b}(\Gamma)(-\cD_K)
=M_{2b}(\Gamma)(-2\cD_K)$.  For $F\in M_{2b}(\Gamma)(-2\nu\cD_K)$ we
consider $(F/w^{2\nu})$ as a function on $\cD_K$.  From the definition
of modular form (Definition~\ref{modforms}) it follows that this
function is holomorphic, $\Gamma\cap\Gamma_K$-invariant (see equation
(\ref{GammaK})) and homogeneous of degree $2b+2\nu$.  Thus
$(F/w^{2\nu})|_{\cD_K}\in M_{2(b+\nu)}(\Gamma\cap\Gamma_K)$.  In
Proposition~\ref{branchdiv} we saw that, $\Gamma_K$ contains
$\Tilde{\Orth}^+(K)$ as subgroup of $\Tilde\Orth^+(L)$(with equality in almost all cases), so
we may replace $\Gamma\cap\Gamma_K$ by $\Gamma\cap\Tilde{\Orth}^+(K)$.
In this way we obtain an exact sequence
$$
0\to M_{2b}(\Gamma)(-(2+2\nu)\cD_K)\to M_{2b}(\Gamma)(-2\nu \cD_K)
\to M_{2(b+\nu)}(\Gamma\cap\Tilde{\Orth}^+(K)),
$$
where the last map is given by $F\mapsto F/w^{2\nu}$. This gives the result.
\end{proof}

Now we proceed with the proof of Theorem~\ref{gentype}.

Let $L$ be a lattice of signature $(2,n)$ and $\Gamma<\Tilde\Orth^+(L)$: 
recall that $k$ is even.
According to Proposition~\ref{obstructions} we can find
pluricanonical  differential forms
on $\cF_L(\Gamma)^{\tor}$ if
\begin{equation}\label{obstruction}
C_B(\Gamma)=\dim M_{k(n-a)}(\Gamma)-\sum_K \dim B(K)>0,
\end{equation}
where summation is taken over all irreducible components of the
ramification divisor (see the remark at the end of \S 2).  It now
remains to estimate the dimension of $B(K)$ for each of the finitely
many components of the ramification locus in the cases we are interested in,
namely $\Gamma=\Orth^+(II_{2,8m+2})$ and $\Gamma=\Tilde\Orth^+(\Ldm)$. 

According to the Hirzebruch-Mumford proportionality principle
$$
\dim M_k(\Gamma)=\frac{2}{n!}\vol_{HM}(\Gamma)k^n+O(k^{n-1}).
$$
The exact formula for the Hirzebruch-Mumford volume $\vol_{HM}$ for
any indefinite orthogonal group was obtained in~\cite{GHS1}. It
depends mainly on the determinant and on the local densities of the
lattice $L$.  Here we simply quote the estimates of the dimensions of
certain spaces of cusp forms.

The case of $II_{2,8m+2}$ is easier because the branch divisor has
only one irreducible component defined by any $(-2)$-vector $r$.  
According to Proposition~\ref{orbits} the orthogonal complement $K_r$ 
is $K_{II}^{(m)}$. This lattice differs from the lattice $L_{2}^{(m)}$,
whose Hirzebruch-Mumford volume was calculated in~\cite[\S
3.5]{GHS1}, only by one copy of the hyperbolic plane. Therefore
$$
\vol_{HM}\Tilde{\Orth}^+(L_{2}^{(m)})=(B_{8m+4}/(8m+4))
\vol_{HM}\Tilde{\Orth}^+(K_{II}^{(m)}), 
$$
and hence, for even $k$,
$$
\dim M_{k}(\Tilde{\Orth}^+(K_{II}^{(m)}))=
\frac{2^{1-4m}}{(8m+1)!}
\cdot \frac{B_2\dots B_{8m+2}}{(8m+2)!!}\,k^{8m+1}+O(k^{8m}),
$$
where the $B_i$ are Bernoulli numbers.
Assume that $m\ge 3$. Let us take a cusp form
$$
F\in  S_{4m+12}(\Orth^+(II_{2,8m+2}))
$$
from  Proposition~\ref{cuspforms}.
In this case the dimension of the obstruction space $B$
of Proposition~\ref{obstructions} for  the pluricanonical forms of order
$k=2k_1$ is given by
\begin{eqnarray*}
\lefteqn{\sum_{\nu=0}^{k_1-1}\dim M_{(4m-10)k+2\nu}(\Tilde{\Orth}^+(K_{II}^{(m)}))=}\\
&&\frac{2^{4m+2}}{(8m+2)!}
\cdot \frac{B_2 \dots B_{8m+2}}{(8m+2)!!}
\big(\bigl(1+\frac 1{4m-10}\bigr)^{8m+2}-1\big)
((4m-10)k_1)^{8m+2}\\ &&+O(k^{8m+1})
\end{eqnarray*}
In \cite[\S 3.3]{GHS1} we computed the leading term of the dimension
of the space of modular forms for $\Orth^+(II_{2,8m+2})$.  Comparing
these two we see that the constant $C_B(\Orth^+(II_{2,8m+2}))$ in
the obstruction inequality (\ref{obstruction}) is positive if and only
if
\begin{equation}\label{BII}
\frac{B_{4m+2}}{4m+2}>\bigl(1+\frac 1{4m-10}\bigr)^{8m+2}-1.
\end{equation}
Moreover $\cF_{II}^{(m)}$ is of general type if $C_B(\Orth^+(II_{2,8m+2}))>0$.
From Stirling's formula
\begin{equation}\label{Stirling}
5\sqrt{\pi n}\big(\frac{n}{\pi e})^{2n}
>|B_{2n}|>4\sqrt{\pi n}\big(\frac{n}{\pi e})^{2n}.
\end{equation}
Using this estimate we easily obtain that (\ref{BII}) holds if $m\ge 5$.
Therefore we have proved Theorem~\ref{gentype} for the lattice $II_{2,8m+2}$.

Next we consider the lattice $\Ldm$ of $\Kthree$ type. For
this lattice the branch divisor of $\Fdm$ is calculated in
Proposition~\ref{orbits}. It contains one or two (if $d\equiv 1 \bmod
4$) components defined by $(-2)$-vectors and some number of components
defined by $(-2d)$-vectors. To estimate the
obstruction constant $C_B(\Gamma)$ in (\ref{obstruction}) we use the
dimension formulae for the space of modular forms with respect to the
group $\Tilde\Orth^+(M)$, where $M$ is one of the following lattices
from Proposition~\ref{orbits}: $\Ldm$ (the main group); $\Kdm$ and
$\Ndm$ (the $(-2)$-obstruction); $M_{2,8m+2}$, $\Kmm$ and $T_{2,8m+2}$
(the $(-2d)$-obstruction). The corresponding dimension formulae were
found in [GHS1] (see \S\S 3.5, 3.6.1--3.6.2, 3.3 and 3.4). The branch
divisor of $(-2d)$-type appears only if $d>1$. We note that
\begin{equation}\label{TKcomp}
\vHM(\Tilde\Orth^+(\Tm))>
\vHM(\Tilde\Orth^+(\Kmm)).
\end{equation}
Therefore in order to estimate $C_B(\Gamma)$
we can assume that all $(-2d)$-divisors defined by stably reflective
$(-2d)$-vectors $r$ with $\divv(r)=d$ (see Proposition~\ref{orbits}) are of
the type $\Tm$.

We put $k=2k_1$, $w=n-a$ and $n=8m+3$. For the obstruction
constant in (\ref{obstruction}) we obtain
\begin{equation}\label{obstr_ineq}
C_B(\Tilde\Orth^+(\Ldm))>
\dim M_{2k_1w}(\Tilde\Orth^+(\Ldm))-B_{(-2)}-B_{(-2d)}
\end{equation}
where
$$
B_{(-2)}=\dim B(\Kdm)+\dim B(\Ndm),
$$
$$
B_{(-2d)}= 2^{\rho(d)-1}(\dim B(\Mm)
+2^{h_d}\dim B(\Tm))
$$
and $B(K)$ is the obstruction space from Proposition~\ref{obstructions}.
By $h_d$ we denote the sum $\delta_{0,d(8)}-\delta_{2,d(4)}$, where $d(n)$ is
$d \bmod n$ and $\delta$ is the Kronecker delta
(see Proposition~\ref{orbits} and the remark following it).

For any lattice considered above
\begin{eqnarray*}
B(K)&=&\sum_{\nu=0}^{k_1-1}
\dim M_{2(k_1w+\nu)}(\Tilde\Orth^+(K))\\
&=&\frac{2^{8m+3}}{(8m+3)!}E_w(8m+3)
\vHM(\Tilde\Orth^+(K))(k_1w)^{8m+3}+O(k_1^{8m+2})
\end{eqnarray*}
where $E_w(8m+3)=(1+\frac 1 w )^{8m+3}$.

All terms in (\ref{obstr_ineq}) contain a common factor.  First
\begin{equation}
\dim M_{2k_1w}(\Tilde\Orth^+(\Ldm))= C_{m,d}^{k_1,w}
\left|\frac{B_{8m+4}}{B_{4m+2}}\right| \sqrt{d}
+O(k_1^{8m+2}),
\end{equation}
where
\begin{eqnarray*}
\lefteqn{C_{m,d}^{k_1,w}=}\\
&&\frac{2^{4m+1+\delta_{1,d}}}{(8m+3)!}
\frac{|B_2\dots B_{8m+2}|}{(8m+2)!!}
\frac{|B_{4m+2}|}{4m+2} d^{4m+\frac{3}{2}}
\prod_{p\mid d}(1-p^{-(4m+2)})(k_1w)^{8m+3}.
\end{eqnarray*}
We note that $2^{4m+1}\frac{B_{4m+2}}{4m+2}=
\pi^{-(4m+2)}\Gamma(4m+2)\zeta(4m+2)$.

From \cite[(16)]{GHS1} it follows that
\begin{eqnarray*}
\lefteqn{\vHM(\Tilde\Orth^+(\Kdm))=}\\
&&2^{\delta_{1,d}-\delta_{4,d(8)}}
\frac{B_2\dots B_{8m+2}}{(8m+2)!!}d^{4m+\frac{3}{2}}
\pi^{-(4m+2)}\Gamma(4m+2)L(4m+2,\left(\frac{4d}{*}\right)).
\end{eqnarray*}
We can use the formula for the volume of $\Ndm$
in the following form:
\begin{eqnarray*}
\lefteqn{\vHM(\Tilde\Orth^+(\Ndm))=}\\
&&2^{1+\delta_{1,d}-(8m+4)}
d^{4m+\frac{3}{2}}
\frac{B_2\dots B_{8m+2}}{(8m+2)!!}\pi^{-(4m+2)}
\Gamma(4m+2)L(4m+2,\left(\frac{d}{*}\right))
\end{eqnarray*}
(see \cite[3.6.2]{GHS1}). It follows that
$$
B_{(-2)}=C_{m,d}^{k_1,w}E_w(8m+3)
(2^{8m+3-\delta_{4,d(8)}} P_K(4m+2)+P_N(4m+2))
+O(k_1^{8m+2})
$$
where
$$
P_K(n)=(1-2^{-n})^{\delta_{0,d(2)}}
\frac{L(n,\left(\frac{4d}{*}\right))}
{L(n,\chi_{0,4d})}\prod_{p\mid d} \frac{1-p^{-n}}{1+p^{-n}}
$$
and
$$
P_N(n)=\frac{L(n,\left(\frac{d}{*}\right))}
{L(n,\chi_{0,d})}\prod_{p\mid d} \frac{1-p^{-n}}{1+p^{-n}}.
$$
Here $\chi_{0,f}$ denotes the principal Dirichlet
character modulo $f$.

We note that $|P_K(n)|<1$ and $|P_N(n)|<1$ for any
$d$. We conclude that
$$
B_{(-2)}<C_{m,d}^{k_1,w} E_w(8m+3) b_{(-2)}
$$
where $b_{(-2)}=2^{8m+3}-1$.

The $(-2d)$-contribution is calculated according to
\cite[3.3--3.4]{GHS1}. We note that $\Tilde\Orth^+(\Tm)$
is a subgroup of $\Tilde\Orth^+(\Mm)$. We obtain
$$
B_{(-2d)}\le C_{m,d}^{k_1,w} E_w(8m+3)b_{(-2d)}
$$
where for $d>2$
$$
b_{(-2d)}=\frac{2^{\rho(d)}}{d}
\left(\frac{4}{d}\right)^{4m+\frac{1}{2}}
4(2^{h_d}(1+2^{-(4m+2)}-2^{-(8m+3)})+2^{-(8m+3)}).
$$
As a result we see that
that the obstruction constant $C_B(\Tilde\Orth^+(\Ldm))$ is positive
if
$$
\beta^{(w)}_{m,d}=\left|\frac{B_{4m+2}}{B_{8m+4}}\right|
E_w(8m+3)(b_{(-2)}+b_{(-2d)})<\sqrt{d}.
$$
Using (\ref{Stirling}) we get
$$
\left|\frac{B_{4m+2}}{B_{8m+4}}\right|
<\frac{5}{4\sqrt{2}}\left(\frac{\pi e}{2m+1}\right)^{4m+2}
\frac{1}{2^{8m+4}}.
$$
For $m\ge 5 $ we choose a cusp form $F_a$ of weight
$a=4m+10$, i.e. we take $w=4m-7$ in Proposition~\ref{obstructions}.
Such a cusp form exists for all $d\ge 1$ by Proposition~\ref{cuspforms}.
Using the fact that $\beta_{(-2d)}\le
\beta_{(-4)}$ for any $d\ge 2$ and the value $b_{(-4)}=2^{4m+\frac{5}{2}}$,
we see that
$$
\beta^{(4m-7)}_{m,d}<(1+\frac{1}{4m-7})^{8m+3}
\frac{5}{8\sqrt{2}}\left(\frac{\pi e}{2m+1}\right)^{4m+2}
\frac{2^{8m+3}+2^{4m+\frac{5}{2}}+1}{2^{8m+3}},
$$
which is smaller than $1$ if $m\ge 5$.
This proves Theorem~\ref{gentype} for $m\ge 5$.

For $m=4$ there exists a cusp form $F_a$ of weight $4m+6$ if $d\neq
1,2, 4$, {i.e.} we take $w=4m-3$. To see that
$\beta^{(13)}_{4,d}<\sqrt{d}$ we need check this only for $d=3$
because $b_{(-2d)}<b_{(-6)}$ for $d>3$.  One can do it by direct
calculation.

For $m\le 3$ we choose $F_a$ of weight $4m+2$, {i.e.}  we take
$w=4m+1$. Such a cusp form exists if $d>180$ according
to Proposition~\ref{cuspforms}. For such $d$ we see that $\beta_{(-2d)}<1$.
Then the obstruction
constant $C_B(\Tilde\Orth^+(\Ldm))$ is positive if
$$
\left|\frac{B_{4m+2}}{B_{8m+2}}\right|
(1+\frac{1}{4m+1})^{8m+3} (2^{8m+3}+2)<\sqrt{d}.
$$
This inequality gives us the bound on $d$ in Theorem~\ref{gentype}. 

This completes the proof of Theorem~\ref{gentype}. 
\bigskip

In the proof of Theorem~\ref{gentype} above we have seen that the $(-2)$-part
of the branch divisor forms the most important reflective obstruction to the extension
of the $\Tilde\Orth^+(\Ldm)$-invariant differential forms to a smooth
compact model of $\Fdm$. Let us consider the double covering $\SFdm$ of
$\Fdm$ for $d>1$ determined by the special orthogonal group:
$$
\SFdm=\SOrthdm\setminus \Ddm\to \Fdm.
$$
Here the branch divisor does not contain
the $(-2)$-part. Theorem~\ref{sgentype} below shows that there are
only five exceptional varieties $\cS\cF_{2d}^{(m)}$
with $m>0$ and $d>1$ that are possibly not of general type.

The variety $\cS\cF_{2d}^{(2)}$ can be interpreted as the moduli space
of $\Kthree$ surfaces of degree $2d$ with spin structure:
see \cite[\S 5]{GHS2}.
The three-fold $\cS\cF_{2d}^{(0)}$ is the moduli space
of $(1,t)$-polarised abelian surfaces.

\begin{theorem}\label{sgentype}
The variety $\SFdm$ is of general type for any $d>1$ if $m\ge 3$.
If $m=2$ then $\cS\cF_{2d}^{(2)}$ is of general type if
$d\ge 3$. If $m=1$ then $\cS\cF_{2d}^{(1)}$ is of general type if
$d=5$ or $d\ge 7$.
\end{theorem}
\begin{proof}
The case $m=2$ is \cite[Theorem 5.1]{GHS2}, and the result for 
$m\ge 5$ is immediate from Theorem~\ref{gentype}. For $m=1$, $3$ and $4$ we can prove more
than what follows from Theorem~\ref{gentype}.

The branch divisor of $\SFdm$ is defined by the reflections
in vectors $r\in \Ldm$ such  that $-\sigma_{r}\in \SOrthdm$,
because the rank of $\Ldm$ is odd. Therefore $r^2=-2d$, by
Proposition~\ref{rsigma}.

If $F\in M_{2k+1}(\SOrthdm)$ is a modular form (note that the
character $\det$ is trivial), $d>1$
and $z\in \cD^{\bullet}_\Ldm$ is such that $(z,r)=0$, then
$$
F(z)=F(-\sigma_r(z))=F(-z)=(-1)^{2k+1}F(z).
$$
Therefore any modular form of odd weight for $\SOrthdm$ vanishes
on the branch divisor.

To apply the low-weight cusp form trick
used in the proof of Theorem~\ref{gentype}
one needs a cusp form of weight smaller than $\dim \SFdm=8m+3$.
By Proposition~\ref{cuspforms} there exists a cusp form
$F_{11+4m}\in S_{11+4m}(\SOrthdm)$.
For $m\geq 3$ we have that $11+4m<8m+3$. Therefore
the differential forms $F_{11+4m}^kF_{(4m-8)k}(dZ)^k$,
for arbitrary  $F_{(4m-8)k}\in M_{(4m-8)k}(\SOrthtd^+(L_{2d}^{(m)}))$,
extend to the toroidal compactification
of $\SFdm$ constructed in Theorem~\ref{torcomp}.
This proves the cases $m\ge 3$ of the theorem.

For the case $m=1$ we use a cusp form of weight $9$ with
respect to $\SOrthtd^+(L_{2d}^{(1)})$ constructed in Proposition~\ref{cuspforms}.
\end{proof}

We can obtain some information also for some of the remaining cases.
\begin{proposition}\label{SF8and12}
The spaces $\cS\cF_{8}^{(1)}$ and  $\cS\cF_{12}^{(1)}$ have
non-negative Kodaira dimension.
\end{proposition}
\begin{proof}
By Proposition~\ref{cuspforms} there are cusp forms of weight $11$ for $\SOrthtd^+(L_{8}^{(1)})$
and  $\SOrthtd^+(L_{12}^{(1)})$. The weight of these forms is equal to the dimension. 
By the well known criterion of Freitag these cusp forms determine
canonical differential forms on the $11$-dimensional varieties $\cS\cF_8^{(1)}$ and $\cS\cF_{12}^{(1)}$.
\end{proof}
It may be that these varieties have intermediate Kodaira dimension.

In \cite{GHS2} we used pull-backs of the Borcherds modular
form $\Phi_{12}$ on $\cD_{II_{2,26}}$ to show that many moduli spaces
of $\Kthree$ surfaces are of general type. We can also use Borcherds
products to prove results in the opposite direction.
\begin{theorem}\label{FII1-2}
The Kodaira dimension of $\cF_{II}^{(m)}$ 
is $-\infty$ for $m=0$, $1$ and $2$.
\end{theorem}
\begin{proof}
  For $m=0$ we can see immediately that the quotient is rational: a
  straightforward calculation gives that $\cF_{II}^{(0)}= \Gamma
  \backslash \HH_1 \times \HH_1$ where $\HH_1$ is the usual upper half
  plane and $\Gamma$ is the degree $2$ extension of $\SL(2,\ZZ) \times
  \SL(2,\ZZ)$ by the involution which interchanges the two
  factors. Compactifying this, we obtain the projective plane $\PP_2$.

  For $m=1$, $2$ we argue differently. 
  There are modular forms similar to $\Phi_{12}$ for the even
  unimodular lattices $II_{2,10}$ and $II_{2,18}$.  They are Borcherds
  products $\Phi_{252}$ and $\Phi_{127}$ of weights $252$ and $127$
  respectively, defined by the automorphic functions
$$
\Delta(\tau)^{-1}(\tau)E_4(\tau)^2=q^{-1}+504+q(\dots)
$$
and
$$
\Delta(\tau)^{-1}(\tau)E_4(\tau)=q^{-1}+254+q(\dots),
$$
where $q=\exp(2\pi i \tau)$ and $\Delta(\tau)$ and $E_4(\tau)$ are
the Ramanujan delta function and the Eisenstein series of weight $4$
(see \cite{B1}).  The divisors of $\Phi_{252}$ and $\Phi_{127}$
coincide with the branch divisors of $\cF_{II}^{(1)}$ and
$\cF_{II}^{(2)}$ defined by the $(-2)$-vectors. Moreover $\Phi_{252}$
and $\Phi_{127}$ each vanishes with order one along the respective
divisor. Therefore if $F_{10k}(dZ)^k$ (or $F_{18k}(dZ)^k$) defines a
pluricanonical differential form on a smooth model of a toroidal compactification
of $\cF_{II}^{(1)}$ or $\cF_{II}^{(2)}$, then $F_{10k}$ (or
$F_{18k}$) is divisible by $\Phi_{252}^k$ (or $\Phi_{127}^k$), 
since $F_{10k}$ or $F_{18k}$ must vanish to order at
least $k$ along the branch divisor. This is not possible, because the
quotient would be a holomorphic modular form of negative weight.
\end{proof}

\bibliographystyle{alpha}

\bigskip
\noindent
V.A.~Gritsenko\\
Universit\'e Lille 1\\
Laboratoire Paul Painlev\'e\\
F-59655 Villeneuve d'Ascq, Cedex\\
France\\
{\tt valery.gritsenko@math.univ-lille1.fr}
\bigskip

\noindent
K.~Hulek\\
Institut f\"ur Algebraische Geometrie\\
Leibniz Universit\"at Hannover\\
D-30060 Hannover\\
Germany\\
{\tt hulek@math.uni-hannover.de}
\bigskip

\noindent
G.K.~Sankaran\\
Department of Mathematical Sciences\\
University of Bath\\
Bath BA2 7AY\\
England\\
{\tt gks@maths.bath.ac.uk}

\end{document}